\documentclass[%
 aip,
 amsmath,amssymb,
preprint,%
]{revtex4-1}

\usepackage{graphicx}
\usepackage{dcolumn}
\usepackage{bm}
\usepackage{subfig}
\usepackage[utf8]{inputenc}
\usepackage[T1]{fontenc}
\usepackage{mathptmx}
\usepackage{bigints}
\usepackage{amsthm}
\usepackage{etoolbox}
\newtheorem{definition}{Definition}
\newtheorem{lemma}{Lemma}
\newtheorem{theorem}{Theorem}

\newcommand{\eye}{\mathbb{I}}
\makeatletter
\def\@email#1#2{%
 \endgroup
 \patchcmd{\titleblock@produce}
  {\frontmatter@RRAPformat}
  {\frontmatter@RRAPformat{\produce@RRAP{*#1\href{mailto:#2}{#2}}}\frontmatter@RRAPformat}
  {}{}
}%
\makeatother
\begin{document}
\title[]{Phase transition in a kinetic mean-field game model of inertial self-propelled agents}
\author{Piyush Grover}
\email{piyush.grover@unl.edu}
 \affiliation{University of Nebraska-Lincoln}
\author{Mandy Huo}%
\affiliation{ 
California Institute of Technology
}%

\date{\today}

\begin{abstract}
The framework of Mean-field Games (MFGs) is used for modelling the collective dynamics of large populations of non-cooperative decision-making agents. We formulate and analyze a kinetic MFG model for an interacting system of non-cooperative motile agents with inertial dynamics and finite-range interactions, where each agent is minimizing a biologically inspired cost function. By analyzing the associated coupled forward-backward in time system of nonlinear Fokker-Planck and Hamilton-Jacobi-Bellman equations, we obtain conditions for closed-loop linear stability of the spatially homogeneous MFG equilibrium that corresponds to an ordered state with non-zero mean speed. Using a combination of analysis and numerical simulations, we show that when energetic cost of control is reduced below a critical value, this equilibrium loses stability, and the system transitions to a travelling wave solution. Our work provides a game-theoretic perspective to the problem of collective motion in non-equilibrium biological and bio-inspired systems.
\end{abstract}

\maketitle


\begin{quotation}
How can we understand phase transitions in collective motion of a population of motile agents making individually optimal decisions over time ? Mean field game (MFG) theory is a mathematical framework for modelling such systems. We study a system consisting of a large number of agents that optimally modulate their acceleration to minimize the long term average of a weighted sum of two costs. The first is the total energy use, and the second penalizes the mismatch between the local mean speed and a preferred speed. We show that the system transitions from a spatially homogeneous ordered state to a travelling wave as the energetic cost is reduced below a critical value.

\end{quotation}

\section{Introduction}
This paper is concerned with phase transitions in the collective motion of optimal decision-making non-cooperative motile agents. Collective motion and synchronization continues to be actively studied in many scientific domains\cite{vicsek1995novel,levine2000self,strogatz2000kuramoto,cucker2007emergent,gompper20202020,kassabov2021sufficiently}. In purely physical systems, the interactions between particles are given by physical laws, and systems are studied using agent-based (or microscopic), kinetic and hydrodynamics models that incorporate those laws. On the other hand, models for systems involving decision-making \emph{agents}, e.g.,  bird flocking, traffic, human crowds and robot swarms, often incorporate other types of interactions such as collision avoidance, alignment, and cohesion towards the mean. Most commonly, these latter types of interactions are captured in phenomenological models using physical analogies.

A different, `inverse modelling' approach \cite{ding2022mean} is to stipulate that the collective behavior of a population of decision-making agents is a solution to a collective optimization or optimal control problem. When considering large number of non-cooperative agents making sequential decisions, the framework of Mean-Field Games (MFGs) \cite{caines2021mean} is appropriate. In a MFG system, the collective behavior is the result of each agent solving an optimal control problem that depends on its own state and control as well as the collective state \cite{bensoussan2013mean,Huang2007}. MFGs formulated in continuous state space and time are described by coupled set of forward-backward in time nonlinear partial differential equations (PDEs). While standard kinetic or hydrodynamic equations used for modelling collective behavior are initial value problems (IVP or evolution PDEs) , the MFG systems have a forward-backward in time structure, and hence consist of boundary value problem (BVP in time PDEs) .


In evolution PDEs, phase transitions \cite{solon2015pattern,strogatz2000kuramoto,orosz2005bifurcations} in collective behavior are studied via stability and bifurcation analysis of steady or time-periodic states. In recent years, similar questions are being explored in the context of MFGs, where bifurcations can be interpreted as phase transitions in the collective behavior of optimal decision-making agents as problem parameters such as the unit cost of control are varied. While uniqueness of solutions in MFGs is guaranteed under certain monotonicity conditions \cite{bardi2019non}, there is an increasing interest \cite{nourian2011mean,grover2018mean,cirant2019existence,cirant2018variational,ullmo2019quadratic,salhab2019collective,hongler2020mean,lori2024topological} in studying MFGs that don't satisfy such conditions and hence are expected to have multiple co-existing steady states. 

In this context, starting with the seminal work of Yin et al. \cite{yin2012synchronization}, several authors have studied the phase transition from a `disordered' steady state to an `ordered' travelling wave solution in the MFG problem involving synchronization of oscillators\cite{carmona2023synchronization,cesaroni2024stationary}. The system involves globally cost-coupled agents  with first-order dynamics (i.e., without inertia), moving on a circular domain. Hence, these MFG models can be understood as `inverse models' for the Kuramoto dynamics \cite{kuramoto1975international}. Other works \cite{cirant2019existence,cirant2018variational} have proved the existence of time-periodic solutions in certain classes of MFGs set in periodic spatial domains. 

In this paper, we study phase transitions in a MFG `inverse model' for the kinetic Czir\'ok model \cite{czirok1999collective}, where unlike all the works discussed above, the population consists of locally interacting inertial agents, i.e., the agent dynamics are \emph{second-order in time}, and the agents control their acceleration. As a result, the ordered states of the Czir\'ok model include spatially homogeneous equilibria and spatially inhomogeneous travelling waves. We show that the MFG model we derive inherits these features. Specifically, we provide a stability analysis of the spatially homogeneous equilibria of the kinetic MFG \cite{degond2014large}. Further, we numerically show the phase transition in the MFG to travelling waves as the equilibrium states lose stability when the unit cost of control is reduced below a critical value. Our analysis combines techniques use to study generalized Kuramoto models with inertia \cite{barre2016bifurcations,acebron2000synchronization}, and invariant subspace methods employed in the study of Riccati equations \cite{laub1991invariant,abou2012matrix} as well as linear-quadratic (LQ) MFGs \cite{huang2019linear}. 

\section{Czir\'ok model revisited}\label{sec:garnier}
In this section, we briefly review the Czir\'ok \cite{garnier2019mean} model, and provide a new stability analysis which sets the stage for the stability analysis of the MFG model in the later sections. The model is a biologically inspired interacting particle system of $n$ agents, where each agent is moving on a 1D periodic domain $[0,l]$. The second-order dynamics of the $i$th agent are given by the stochastic differential equation (SDE):
\begin{align}
\textrm{d} x_i &= u_i \textrm{d}t, \\
\textrm{d}u_i &= [G(\langle u \rangle_i) - u_i]\textrm{d}t + \sigma \textrm{d}w_i(t),
\end{align}
    
where $\langle u \rangle_i = \dfrac{1}{n} \displaystyle\sum^n_{j = 1} u_j \phi(\| x_j - x_i\|)$, and 
    $\| x_j - x_i\| = \min (|x_j - x_i|, l - |x_j - x_i|)$ is the distance on the torus  between the positions of the $i^{\textrm{th}}$ and $j$th agents. The $w_i$ are independent Brownian motions, and $\sigma$ is noise intensity.
The finite-range interaction kernel $\phi(x)=\dfrac{l}{2}\times\mathbf{1}_{[0,1]}$ with $\frac{1}{l}\bigintssss^l_0 \phi (\|x\|) \textrm{d}x = 1$,  and $G(u)=\dfrac{h+1}{5}u-\dfrac{h}{125}u^3$. Here,  $G(u)-u$ is the gradient of a bistable potential whose critical points are the preferred mean speeds. In the $n\rightarrow\infty$ limit, the nonlinear PDE describing the evolution of density $\rho(t,x,u)$ is given by the Fokker-Planck equation
\begin{eqnarray}
\frac{\partial{\rho(t,x,u)}}{\partial t}=&&\nonumber\-u\frac{\partial{\rho(t,x,u)}}{\partial x}\nonumber -\frac{\partial}{\partial u}\left( \left[G\left(\int_0^L\int_{-\infty}^{\infty} u'\phi(\|x'\|)\rho(t,x-x',u')du'dx'\right)-u\right]\rho(t,x,u)\right)\\ && + \frac{1}{2}\sigma^2\frac{\partial^2\rho(t,x,u)}{\partial u^2}.\label{eq:PDE_for}
\end{eqnarray} Note that this equation is degenerate since there is no diffusion in the spatial ($x$) variable. The system has spatially homogeneous equilibrium states $\rho_{\xi}(u)$ of the form 
$\rho_{\xi}(u)=\dfrac{1}{l}F_{\xi}(u)$, where $F_{\xi}(u)=\dfrac{1}{\sqrt{\pi\sigma^2}}e^{-(u-\xi)^2/{\sigma^2}}$, and $\xi$ is the mean speed of the population satisfying the fixed point equation $\xi=G(\xi)$. This fixed point equation has a unique `disordered' equilibrium with $\xi=0$ for $h<4$, and two additional ordered equilibria with $\pm\xi\neq 0$ for $h>4$.

\subsection{Linearization and the operator eigenvalue equation}\label{sec:for_lin}

 We focus our stability analysis on the ordered ($\xi\neq 0$) spatially homogeneous equilibria which exist when $h>4$. It is known that for $h>4$, the disordered equilibrium is unstable. We linearize Eq. \ref{eq:PDE_for} around $\rho_{\xi}(u)$, and note that $F'_{\xi}(u)=-\dfrac{2(u-\xi)}{\sigma^2}F_\xi(u)$. We take the perturbed state to be of the form
\begin{align}
\rho(t,x,u)=\rho_{\xi}(u)+\epsilon\sqrt{\rho_{\xi}(u)}\bar{\rho}(t,x,u),\label{eq:for_perdef}
\end{align} and 
expand the perturbation $\bar{\rho}$ using spatial Fourier decomposition:
\begin{align}
\bar\rho(t,x,u)=\sum_{k=0,\pm 1,\pm 2,...} \hat\rho_k(t,u)e^{i2\pi kx/l}.\label{eq:for_perfourier}
\end{align}
We will work in the space of square-integrable complex-valued functions on the real line, i.e, we will assume $\hat{\rho}(t,.)\in L^2(\mathbb{R},du)$, with inner product $\langle a,b\rangle=\bigintssss_{\infty}^{\infty}\bar{a}(u)b(u)du$. The resulting linearized evolution equation for the $k$th Fourier mode is
\begin{align}
    \frac{\partial{\hat\rho_k(t,u)}}{\partial t}=L^k\hat{\rho}_k=(L^k_{loc,1}+L^k_{nonloc})\hat{\rho}_k,
\end{align}
where
\begin{align}
    L^k_{loc,1}[f](u)=\left[\frac{-i2\pi ku}{l}-\frac{[(u-\xi)^2-\sigma^2]}{2\sigma^2}\right]f(u)+\frac{1}{2}\sigma^2\frac{\partial^2 f(u)}{\partial u^2},\label{op_Garnier_loc}\\ 
    L^k_{nonloc}[f](u)=\frac{2(u-\xi)}{\sigma^2}\sqrt{F_{\xi}(u)}G'(\xi)\phi_k\left(\int u'\sqrt{F_{\xi}(u')}f(u')du'\right).\label{op_Garnier_nonloc}\end{align}

The eigenfunctions and eigenvalues of the operator $L^k_{loc,1}$ can be computed explicitly as follows. We perform a change of variables \cite{acebron2000synchronization} $u=\mathfrak{g}(v)=a_1v+a_2$, where $a_1=\sigma/\sqrt{2}$ , and $a_2=\xi-2ik\pi\sigma^2/l$. With this substitution, the eigenvalue equation for $L^k_{loc,1}$ reduces to \begin{align}
\alpha_{k}\hat\rho_{loc,k}=\left[\left(\frac{1}{2}-\frac{v^2}{4}+c_2(k)\right)\hat\rho_{loc,k}+\frac{\partial^2 \hat\rho_{loc,k}}{\partial v^2}\right]. \label{eq:eig_Garnier_locv}
    \end{align}
It is known \cite{AbramowitzStegun1964} that the eigenvalue Eq. \ref{eq:eig_Garnier_locv} has solutions 
$\alpha_{k,p}=-p+c_2(k), p=0,1,2,\dots$. Here $c_2(k)=-\dfrac{2\pi^2k^2\sigma^2}{l^2}-\dfrac{2ik\pi\xi}{l}.$ The corresponding eigenfunctions are parabolic cylinder functions $D_p(v)=2^{-\frac{p}{2}}e^{\frac{-v^2}{4}}H_p(\dfrac{v}{\sqrt{2}})=e^{\frac{-v^2}{4}}\tilde{H}_p(v)$. Here $H_p$ are the physicist's Hermite polynomials and $\tilde{H}_p$ are the probabilist's Hermite polynomials.
Hence, operator $L^k_{loc,1}$ also has the same eigenvalues $\alpha_{k,p}$. Its eigenfunctions are $\eta_{k,p}(u)=z(p)D_p(\mathfrak{g}^{-1}(u))=z(p)D_p(\dfrac{\sqrt{2}}{\sigma}(u-\xi+\dfrac{2ik\pi\sigma^2}{l}))$. Here $z(p)=\sqrt{\dfrac{1}{\sqrt{2\pi}a_1 p!}}$ is a normalizing factor.
The adjoint $(L^k_{loc,1})^*$ has eigenvalues $\bar{\alpha}_{k,p}$, and eigenfunctions $\psi_{k,p}=\bar{\eta}_{k,p}$. Note that due to our chosen normalization, $\langle \eta_{k,p},\psi_{k,q}\rangle=\delta_{pq}$. Hence the set $\{\eta_{k,p},\psi_{k,q}\}_{p=0,1,\dots,q=0,1,\dots}$ form a biorthogonal eigenfunction basis set of $L^2[\mathbb{R},du]$.

With the eigenfunctions and eigenvalues of the local operator at hand, we can derive the characteristic equation of full operator $L^k$ as follows. Let ($\lambda,\hat{\rho}$) be an eigenvalue-eigenfunction pair of $L^k$. The eigenvalue equation is
    \begin{align}
    \lambda \hat{\rho}(u)=L^k_{loc,1}[\hat{\rho}](u)+g_0(u)\langle \bar{\hat{\rho}},s\rangle, \label{eq:foreig}
    \end{align}
where we have defined functions $g_0(u)\triangleq\dfrac{2(u-\xi)}{\sigma^2}\sqrt{F_{\xi}(u)}G'(\xi)\phi_k$, and $s(u)\triangleq u\sqrt{F_{\xi}(u)}$.

\begin{figure}[h!]
\subfloat[]{
\includegraphics[width=0.245\textwidth]{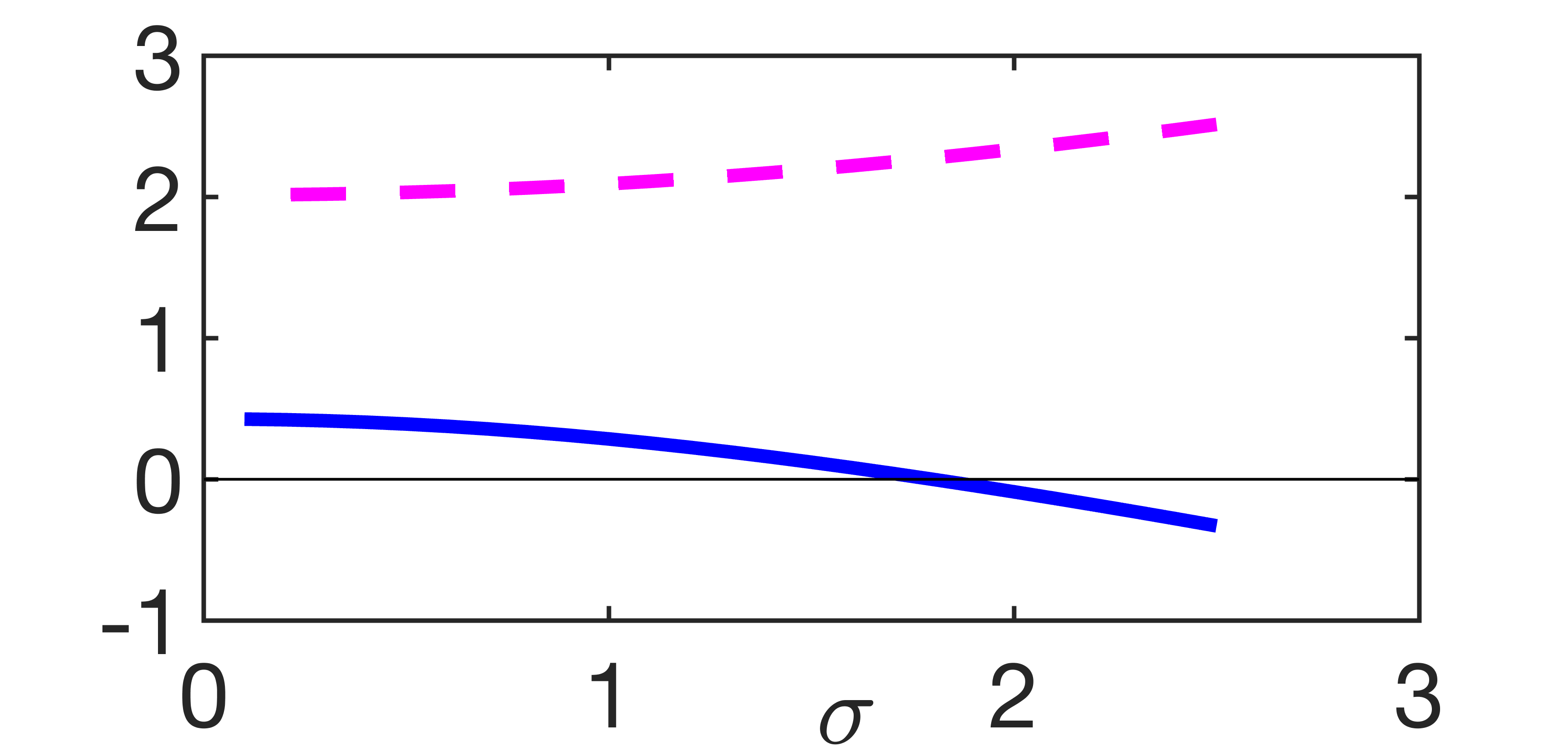}}
\subfloat[]{
\includegraphics[width=0.245\textwidth]{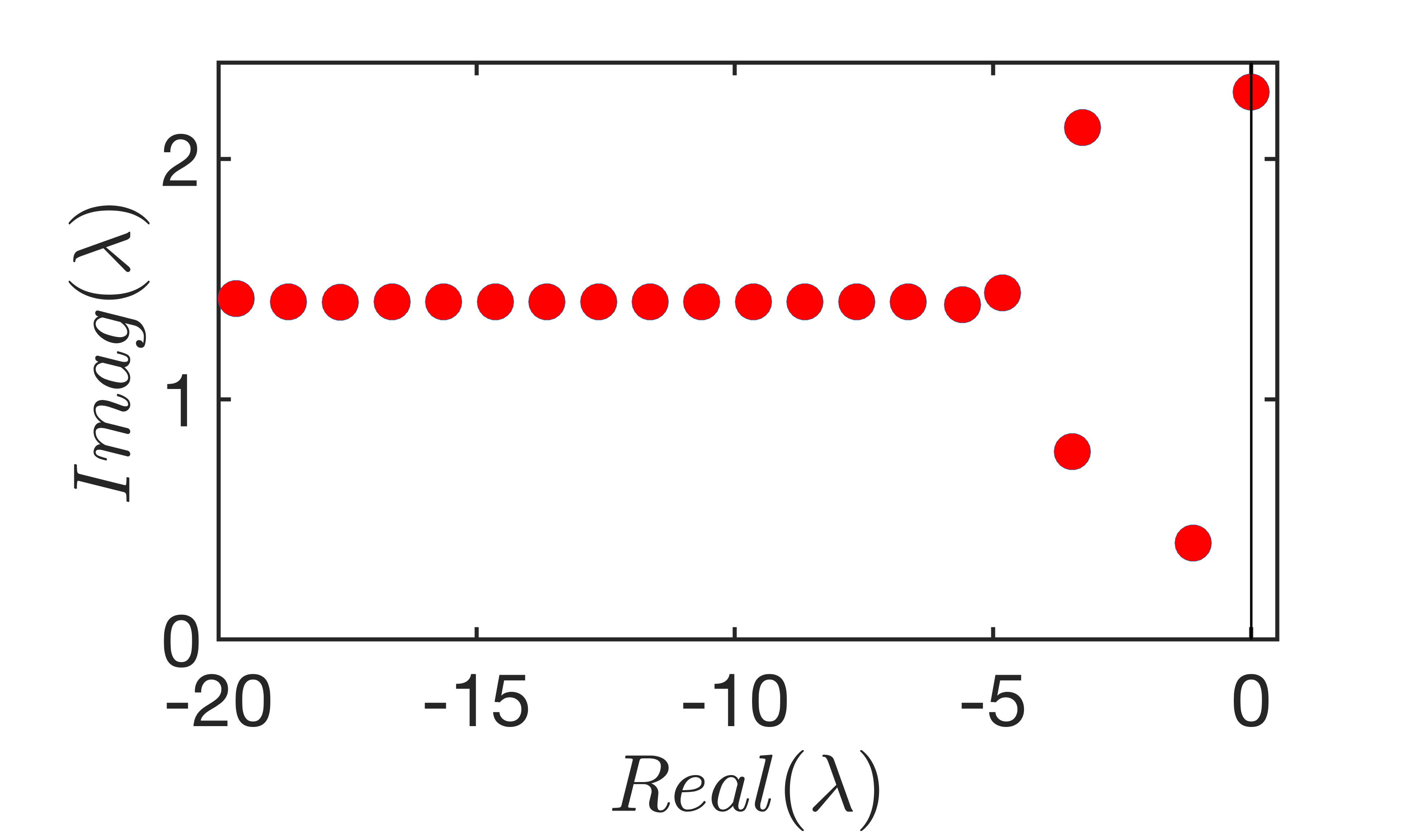}}
\caption{\footnotesize{ (a) The real (solid) and imaginary (dashed) components of the eigenvalue of $L_{k}$ (for $h=5,k=1$) closest to the imaginary axis, as a function of noise intensity $\sigma$. The system is unstable for $\sigma\leq \sigma_c=1.8$ as the real part is positive in that range. (b) The spectrum of $L_{k}$ at the threshold of stability, $\sigma=\sigma_c$. }}
\label{fig:for_eval}
\end{figure}

Let $R^k_{1,\lambda}=(L^k_{loc,1}-\lambda)^{-1}$ be the resolvent of $L^k_{loc,1}$. Then, the action of $R^k_{1,\lambda}$ on an arbitrary complex-valued function f(u) is 
$R^k_{1,\lambda}[f](u)=\displaystyle\sum_{p=0}^{\infty}\dfrac{\langle \psi_{k,p}, f\rangle }{\alpha_{k,p}-\lambda}\eta_{k,p}(u)$. Using this relation, Eq. \ref{eq:foreig} can be rewritten as $\hat{\rho}= -\displaystyle\sum_{p=0}^{\infty} \dfrac{ \langle \psi_{k,p}, g_0 \rangle \langle\bar{\hat{\rho}},s\rangle }{\alpha_{k,p}-\lambda} \eta_{k,p}(u)$.
Taking the inner product of the conjugate of this expression with $s$, cancelling common terms on both sides, and another conjugate operation yields the characteristic equation for $L^k_{loc,1}$:

\begin{align}\label{eq:for2}
1=\sum_{p=0}^{\infty}\frac{\langle \psi_{k,p}, g_0\rangle \langle \bar{\eta}_{k,p},s\rangle}{\lambda-\alpha_{k,p}}=\sum_{p=0}^{\infty}\frac{\langle \bar{\eta}_{k,p}, g_0\rangle \langle \bar{\eta}_{k,p},s\rangle}{\lambda-\alpha_{k,p}}.
\end{align}
To verify our analysis, we solve Eq. \ref{eq:for2} for $l=10, h=5$, and $0.2\leq\sigma\leq 2.5$ using 20 eigenfunctions for each Fourier mode. As shown in Fig. \ref{fig:for_eval}, the $k=1$ Fourier mode loses stability at $\sigma=\sigma_c=1.8$, which matches the prediction in Garnier et al. \cite{garnier2019mean}.

\begin{figure}[h!]
\includegraphics[width=0.45\textwidth]{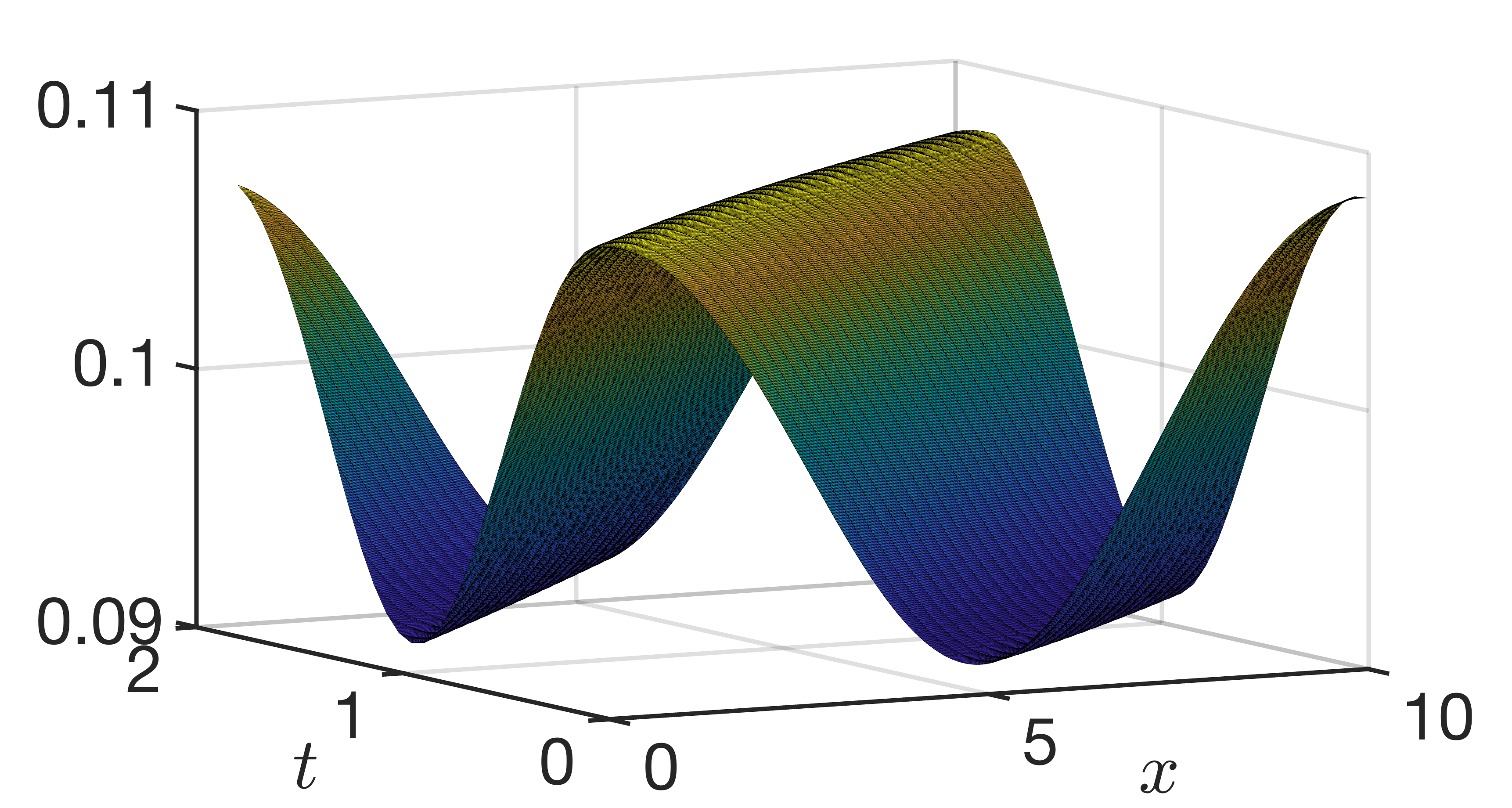}
\caption{\footnotesize{The marginal density $\protect \bigintssss_{\protect \mathbb{R}} \rho(t,x,v)dv$ of a travelling wave solution of the Czir\'ok model Eq. \ref{eq:PDE_for} for $h=5, \sigma=0.8<\sigma_c$. This solution is the steady state reached upon perturbing the unstable spatially homogeneous equilibrium $\protect \rho_{\protect \xi}(u)$, where $\protect \bigintssss_{\protect \mathbb{R}} \rho_{\xi}(v)dv=\protect \dfrac{1}{l}=0.1.$  } }
\label{fig:for_tw}
\end{figure}

\paragraph*{Transition to travelling wave:} When a spatially homogeneous equilibrium $\rho_{\xi}(u)$ ($\xi\neq 0$) loses stability, the system converges to a stable travelling wave \cite{garnier2019mean}, i.e., a time-varying spatially non-homogeneous solution of the form $\rho(t,x,v)=\tilde{\rho}(x-\omega t,v)$. In Fig. \ref{fig:for_tw}, we show one such solution obtained by numerically solving the IVP Eqs. \ref{eq:PDE_for} in the open-source software Dedalus \cite{burns2020dedalus}.

\section{A Mean-field game Czir\'ok model}
To construct the kinetic MFG model, we consider the interacting particle system of $n$ agents, in which  the $i$th agent has the dynamics of the form
    \begin{align}
\textrm{d} x_i &= u_i \textrm{d}t, \\
\textrm{d}u_i &= \alpha_i\textrm{d}t + \sigma \textrm{d}w_i(t),
\end{align}
where $(x_i(t),u_i(t))$ are its position and speed, and $\alpha_i(t)$ is its chosen control.
Each agent seeks to minimize the long term average of the sum of two costs. The first cost depends on a weighted mean speed of the population around the agent's position, and the other is a measure of the control effort. The cost function inspired by the interaction term in the Czir\'ok model Eq. (\ref{eq:PDE_for}) is
\begin{equation}
    J = \limsup_{T \rightarrow \infty} \frac{1}{T} \int^T_0 [c( (x_i, u_i);(x,u)_{-i}) + r \alpha_i^2 ]\textrm{d} t,
\end{equation}
    
with cost-coupling function
\begin{equation}
    c( (x_i, u_i);(x,u)_{-i}) = \left([G(\langle u \rangle_i) - u_i]\right)^2,
\end{equation}
where $(x,u)_{-i}=\{(x_1,u_1),(x_2,u_2),\dots,(x_{i-1},u_{i-1}),(x_{i+1},u_{i+1}),\dots,(x_n,u_n)\}$ represents the rest of the population. Here, $r$ is the unit cost of control. In the $n\rightarrow\infty$ population limit, we can approximate the cost-coupling function as \cite{grover2018mean}
\begin{equation}
    c[\rho](t,x, u) =  \left[G \left(\iint u' \phi (\| x'\|) \rho(t,x - x',u')\textrm{d}u'\textrm{d}x'\right) - u\right]^2.
\end{equation}

Using standard techniques \cite{yin2012synchronization,Nourian2013thesis}, we can show that the corresponding MFG system consist of the following nonlinear Fokker-Planck (FP) and Hamilton-Jacobi-Bellman (HJB) equations governing the density $\rho(x,u,t)$ and the relative value function $h(x,u,t)$ respectively: 

\begin{align}
\frac{\partial{\rho}}{\partial t}=-u\frac{\partial{\rho}}{\partial x}+\frac{1}{2r}\frac{\partial}{\partial u}\{\rho\frac{\partial h}{\partial u}\}+\frac{\sigma^2}{2}\frac{\partial^2\rho}{\partial u^2}\label{eq:mfg_fp},\\
\frac{\partial{h}}{\partial t}=\chi-c[\rho] -u\frac{\partial h}{\partial x}+\frac{1}{4r}(\frac{\partial h}{\partial u})^2 -\frac{\sigma^2}{2}\frac{\partial^2 h}{\partial u^2}\label{eq:mfg_hjb},
\end{align}
where  $\chi$ is the minimum average cost. The PDEs can be solved given an initial density $\rho(0.,.)$, the required final condition $\displaystyle\lim_{t\rightarrow\infty}h(t,.,.,)=0$, and appropriate decay conditions in space. The optimal control $\alpha(t,x,u)=-\dfrac{1}{2r}\dfrac{\partial h}{\partial u}$. It is easily verified that the $(\rho_{\xi^*}(u),h_{\xi^*}(u))$ is a stationary solution to the MFG equations, where $\rho_{\xi^*}(u) = \dfrac{1}{l}F_{\xi^*}(u),
F_{\xi^*}(u)=\dfrac{1}{\sqrt{\pi \sqrt{r}\sigma^2}}e^{-(u-\xi^*)^2/{(\sqrt{r}\sigma^2)}}$ and $
h_{\xi^*}(u)=\sqrt{r}(u-\xi^*)^2$. Here  $\xi^*$ is again the mean-speed of the population, satisfying $G(\xi^*) = \xi^*$, and $\chi = \sigma^2\sqrt{r}$. Hence, similar to Czir\'ok model discussed in Sec. \ref{sec:garnier}, the MFG model also possesses three spatially homogeneous equilibria for $h>4$, among which the two equilibria with $\pm\xi^*\neq 0$ are the ordered states.
    

\subsection{Linearization and the operator eigenvalue equation}\label{sec:lin_mfg}
We linearize the MFG Eqs. (\ref{eq:mfg_fp},\ref{eq:mfg_hjb}) around an ordered equilibrium $(\rho_{\xi}(u),h_{\xi}(u))$ with $\xi=
\xi^*\neq 0$. As in section \ref{sec:for_lin}, we restrict the stability analysis to a class of density perturbations that are exponentially decaying in $u$, i.e., we choose $\rho(t,x,u)=\rho_{\xi}(u)+\epsilon\sqrt{\rho_{\xi}(u)}\bar{\rho}(t,x,u)$, where $\bar{\rho}(t,x,.)\in \mathbb{L}^2\{\mathbb{R},du\}$. In case of the value function, the perturbations are restricted to have a bounded exponential growth rate, i.e., we choose $h(t,x,u)=h_{\xi}(u)+\epsilon\dfrac{\bar{h}(t,x,u)}{\sqrt{\rho_{\xi}(u)}}$, where $\bar{h}(t,x,.)\in \mathbb{L}^2\{\mathbb{R},du\}.$ We expand the perturbations using spatial Fourier decompositions, 
$\bar\rho(t,x,u)=\displaystyle\sum_{k=0,\pm 1,\pm 2,...} \hat\rho_k(t,u)e^{i2\pi kx/l}$, and $\bar h(t,x,u)=\displaystyle\sum_{k=0,\pm 1,\pm 2,...} \hat{h}_k(t,u)e^{i2\pi kx/l}$. With these substitutions, the linearized FP equation is
\begin{align}
\frac{\partial{\hat\rho_k(t,u)}}{\partial t}= -\frac{i2\pi ku}{l}\hat\rho_k-\tilde{L}_{loc}\hat\rho_k-\frac{1}{r\sigma^2}\tilde{L}_{loc}\hat h_k,\label{eq:MFG_lin_FP}
\end{align}
where $\tilde{L}_{loc}\triangleq \left[\dfrac{[(u-\xi)^2-\sqrt{r}\sigma^2]}{2r\sigma^2}\right]-\dfrac{\sigma^2}{2}\dfrac{\partial^2 }{\partial u^2}$ is a self-adjoint operator. Hence, $L^k_{loc,1}= -\dfrac{i2\pi ku}{l}-\tilde{L}_{loc}$.

The linearized HJB equation is 
\begin{align}
\frac{\partial{\hat h}_k(t,u)}{\partial t}=-\frac{i2\pi ku}{l}\hat{h}_k(t,u)+\tilde{L}_{loc}\hat h_k(t,u)-\tilde{c}[\hat{\rho}_k],\label{eq:MFG_lin_HJB}
\end{align}

where $\tilde{c}[\hat{\rho}_k]=2(\xi-u)G'(\xi)\phi_k\sqrt{F_{\xi}(u)}\left(\bigintssss u'\sqrt{F_{\xi}(u')}\hat{\rho}_k(t,u')du'\right)$. We define $L^k_{loc,2}\triangleq -\dfrac{i2\pi ku}{l}+\tilde{L}_{loc}=-(L^k_{loc,1})^*$.  The eigenvalues of $L^k_{loc,2}$ are $\beta_{k,p}=-\bar\alpha_{k,p}=\dfrac{p}{\sqrt{r}}-\dfrac{i2\pi k\xi}{l}+\dfrac{2k^2\pi^2\sigma^2r}{l^2}$. The corresponding eigenfunctions are $\psi_{k,p}(u)$. The coupled system of equations for the linearized MFG are
\begin{align}
\frac{\partial{\hat\rho}_k(t,u)}{\partial t}=L^k_{loc,1}{\hat \rho}_k(t,u)-\frac{1}{r\sigma^2}\tilde{L}_{loc}{\hat h}_k(t,u),\label{eq:lin_MFG_sys1}\\
\frac{\partial{\hat h}_k(t,u)}{\partial t}=-\tilde{c}[\hat{\rho}_k]+L^k_{loc,2}{\hat h}_k(t,u). \label{eq:lin_MFG_sys2}
\end{align}

Let $(\{\hat{\rho}_k,\hat{h}_k\},\lambda)$ be a eigenfunction-eigenvalue pair for the linearized MFG system. The abstract operator eigenvalue equation can be written as:
\begin{align}
\lambda \begin{bmatrix}
\hat{\rho}_k\\ \\ 
\hat{h}_k
\end{bmatrix}
    =  \begin{bmatrix}
L^k_{loc,1}& -\frac{1}{r\sigma^2}\tilde{L}_{loc}\\\\
0 & L^k_{loc,2}
\end{bmatrix} 
\begin{bmatrix}
\hat{\rho}_k\\\\
\hat{h}_k
\end{bmatrix} + \begin{bmatrix}
    0\\\\-g(u)\langle \bar{\hat\rho}_k,s\rangle
\end{bmatrix},\label{eq:eigsys_mfg}
\end{align}
where, $g(u)\triangleq2(\xi-u)\sqrt{F_{\xi}(u)}G'(\xi)\phi_k$, and $s(u)=u\sqrt{F_{\xi}(u)}$. Analogous to the development in Sec. \ref{sec:for_lin}, we employ the resolvents $R^k_{1,\lambda}$ of $L^k_{loc,1}$, and $R^k_{2,\lambda}$ of $L^k_{loc,2}$ to obtain the characteristic equation for the MFG eigenvalues. Here, $R^k_{2,\lambda}[f](u)=\displaystyle\sum_{q=0}^{\infty}\dfrac{\langle \eta_q, f\rangle }{\beta_{k,q}-\lambda}\psi_{k,q}(u)$. Using this relation, the second equation of Eqs. \ref{eq:eigsys_mfg} can be written as \\$\hat{h}_k=\displaystyle\sum_{q=0}^{\infty}\dfrac{\langle \eta_q, g\rangle \langle \bar{\hat{\rho}}_k,s\rangle }{\beta_{k,q}-\lambda}\psi_{k,q}(u)$. Substituting this expression of $\hat{h}_k$ in the first equation of Eqs. \ref{eq:eigsys_mfg} yields
\begin{align}
    R^k_{1,\lambda}[\hat{\rho}_k](v)=\dfrac{1}{r\sigma^2}\sum_{q=0}^{\infty}\dfrac{\langle \eta_q, g\rangle \langle \bar{\hat{\rho}}_k,s \rangle }{\beta_{k,q}-\lambda} \tilde{L}_{loc}[\psi_{k,q}](u).
\end{align}

After applying $(R^k_{1,\lambda})^{-1}$ to above equation, we take the inner product of the conjugate of the resulting equation with $s$. Cancelling common factors from both sides, and taking the conjugate again yields the characteristic equation for the MFG eigenvalues:

\begin{align}
  1=\frac{1}{r\sigma^2}\sum_{p=0}^{\infty}\sum_{q=0}^{\infty}\frac{\langle\bar{\eta}_{k,p},s\rangle \langle \bar{\psi}_q, g\rangle  \langle\bar{\eta}_{k,p}, \tilde{L}_{loc}[\psi_{k,q}]\rangle}{(\alpha_{k,p}-\lambda)(\beta_{k,q}-\lambda)}.\label{eq:mfg_eval}
\end{align}

\subsection{Linear stability of the ordered MFG equilibria}

Since a MFG system is a forward-backward system, the stability of an equilibrium (or more general invariant sets such as periodic orbits or travelling waves) needs to be defined in a way that captures the feedback nature of the system. Consider the linear two-point BVP on the time interval $[0,\infty]$ consisting of the linearized MFG Eqs. (\ref{eq:lin_MFG_sys1},\ref{eq:lin_MFG_sys2}) along with a prescribed initial density perturbation $\hat{\rho}_k(0,.)$, and final value function  $\displaystyle\lim_{t\rightarrow\infty}\hat{h}_k(t,.)= 0$ (inherited from the HJB equation).
\begin{definition}\cite{huang2006large,gueant2009reference,yin2012synchronization,grover2018mean} A MFG equilibrium $(\rho_{\xi},h_{\xi})$ of the MFG system (\ref{eq:mfg_fp},\ref{eq:mfg_hjb}) is said to be linearly stable if the above defined BVP has a unique solution $(\hat{\rho}_k(t,.),\hat{h}_k(t,.))$ such that $\displaystyle\lim_{t\rightarrow\infty}\hat{\rho}_k(t,.)\rightarrow 0$. In other words, any initial density perturbation decays to $0$ in the closed loop as the time horizon goes to $\infty$. \label{def:bvpstability}
\end{definition}

To proceed with the stability analysis, we first obtain the explicit form of the linear BVP discussed above, using the eigenfunction expansions $\hat{\rho}_k(t,u)=\displaystyle\sum_p Y_{1,p}(t)\eta_{k,p}(u)$, and $\hat{h}_k(t,u)=\displaystyle\sum_p Y_{2,p}(t)\psi_{k,p}(u)$. The prescribed initial condition on density fixes the value of $Y_1(0)$, and the final condition on the value function becomes $\displaystyle\lim_{t\rightarrow\infty}Y_2(t)=0$. The ODE system obtained by inserting the expansions into Eqs. (\ref{eq:lin_MFG_sys1},\ref{eq:lin_MFG_sys2}) is 
\begin{align}
\begin{bmatrix}
\dot{Y}_1\\
\dot{Y}_2
\end{bmatrix}
    = N\begin{bmatrix}
Y_1\\
Y_2
\end{bmatrix},\label{eq:bvpmat1}
\end{align} where $N \triangleq \begin{bmatrix}
A_1& B_1\\
A_2 & B_2
\end{bmatrix} $. The submatrices are $(A_1)_{q,p}=\alpha_{k,q-1} \delta_{qp}$, $(A_2)_{q,p}=-\langle \eta_{k,q-1},g(u)\rangle \langle s(u),\eta_{k,p-1}\rangle$,  $(B_1)_{q,p}=\dfrac{-1}{r\sigma^2}\langle\psi_{k,q-1},\tilde{L}_{loc}\psi_{k,p-1}\rangle$, and  $(B_2)_{q,p}=\beta_{k,q-1} \delta_{qp}$. Using the properties of operators involved, we make the following observations. First, $A_1$ is diagonal, $Re(\sigma(A_1))<0$, i.e., all eigenvalues of $A_1$ lie in the left half plane, and $B_2=-A_1^*$. Second, $B_1$ is Hermitian, hence its eigenvalues are real, with $Re(\sigma(B_1))\leq 0$. Finally, $A_2$ is a rank-1 non-Hermitian matrix. The characteristic equation of $N=\Delta+\mathcal{P}_1$ is 
\begin{align}
    det(\Delta+\mathcal{P}_1-\lambda\eye)=det(\Delta-\lambda\eye)det(\eye+(\Delta-\lambda\eye)^{-1}\mathcal{P}_1)=0, \label{eq:mfg_eval_mat}
\end{align} where $\Delta=\begin{bmatrix}
A_1& B_1\\
0 & B_2
\end{bmatrix} $ is an upper block triangular matrix with explicitly known eigenvalues, and $\mathcal{P}_1=\begin{bmatrix}
0 & 0\\
A_2& 0
\end{bmatrix}$ is the rank-1 `perturbation'. Note that this is the explicit form of the characteristic Eq. \ref{eq:mfg_eval}.


Next, we recall some results on  Hamiltonian matrices and Riccati theory \cite{abou2012matrix}. A $2m\times 2m$ complex matrix $H$ is called Hamiltonian if $JH+H^*J=0 $, where $J=\begin{bmatrix}
0 & \eye \\
-\eye & 0
\end{bmatrix}$. If $H$ is Hamiltonian, and $\lambda$ is an eigenvalue of $H$, then $-\bar\lambda$ is also an eigenvalue of $H$. Necessarily, every $2m\times 2m$ Hamiltonian matrix is of the form
\begin{align}
H=\begin{bmatrix}
A & B \\
C & D
\end{bmatrix},\label{eq:ham_mat}\end{align}
where $D=-A^*$, and $B$ and $C$ are both Hermitian. Associated with $H$ is a continuous time algebraic (matrix) Riccati equation (CARE): $XA-DX+XBX-C=0$. We say that an $m\times m$ matrix $X$ that solves the CARE is a \emph{stabilizing solution} if $A+BX$ is Hurwitz. We will need the following result:
\begin{theorem}(Chen et al.\cite{chen2018linear})
Let $H$ be a Hamiltonian matrix of Eq. \ref{eq:ham_mat}, with $A$ Hurwitz (i.e, all its eigenvalues of $A$ lie in the left half plane). Suppose $H$ does not have any eigenvalues on the imaginary axis. Then, there exists an orthogonal (`Schur') transformation $V$ such that \begin{align}
V^*HV= \begin{bmatrix}
H_{11} & H_{12} \\
0 & H_{22}
\end{bmatrix},\end{align}
where all eigenvalues of $H_{11}$ are the stable eigenvalues of $H$. If we block partition $V=\begin{bmatrix}
V_{11} & V_{12} \\
V_{21} & V_{22}
\end{bmatrix}$, then $\begin{bmatrix}
V_{11}\\
V_{21}
\end{bmatrix}$ are the $m$ stable Schur vectors corresponding to the stable block $H_{11}$. Furthermore, the matrix $V_{11}$ is invertible, and $X_+=V_{21}V_{11}^{-1}$ is the unique Hermitian stablizing solution of CARE.
\end{theorem}

Note that the matrix $N$ in Eq. \ref{eq:bvpmat1} is not Hamiltonian since $A_2$ is not Hermitian. We consider its `symmetrized' version, $N_s=\begin{bmatrix}
A_1 & B_1 \\
\dfrac{A_2+A^*_2}{2} & B_2
\end{bmatrix}$, which is a Hamiltonian matrix. Using the definitions of the submatrices discussed above, it follows that the characteristic equation of $N_s$ reduces to that of $N$, i.e., Eq. \ref{eq:mfg_eval_mat}. Hence $N$ and $N_s$ are similar matrices, i.e., they have the same eigenvalues. This implies that there exists an invertible $2m\times 2m$ complex matrix $P$ s.t. $N_s=P^{-1}NP$. Note $P$ can be numerically obtained by eigenvalue decomposition of $N$ and $N_s$. Then we have the following result:

\begin{lemma} Suppose $N$ defined in Eq. \ref{eq:bvpmat1} does not have any eigenvalues on the imaginary axis. Then, the BVP defined by Eq. \ref{eq:bvpmat1} along with arbitrary initial condition $Y_1(0)=Y_{10}$, and final condition $\displaystyle\lim_{t\rightarrow\infty}Y_2(t)=0$ has a solution $(Y_1(t),Y_2(t))$ s.t. $\displaystyle\lim_{t\rightarrow\infty}Y_1(t)=0$.

\begin{proof}
Since $N_s$ and $N$ are similar, $N_s$ also doesn't have any eigenvalues on the imaginary axis. Then $N_s$ satisfies the conditions of Theorem 1, since $A_1$ is Hurwitz. Hence there exists a Hermitian stabilizing solution $X_+$ for CARE corresponding to $N_s$. Define $U=\begin{bmatrix}
\eye &0\\
X_+ &\eye
\end{bmatrix}$, such that $U^{-1}=\begin{bmatrix}
\eye& 0\\
-X_+ &\eye
\end{bmatrix}$. Then, $U^{-1}N_sU=\begin{bmatrix}
A_c &B\\
0 &-A_c^*
\end{bmatrix}$, where $A_c=A_1+B_1X_+$ is a stable matrix by Theorem 1, and we have used the fact the $X_+$ solves the CARE. Now define new states $(Z_1,Z_2)$ via $PU\begin{bmatrix}
Z_1\\
Z_2
\end{bmatrix}=\begin{bmatrix}
Y_1\\
Y_2
\end{bmatrix}$. Note that $PU$ is invertible since it is the product of two invertible matrices. The BVP system of Eq. \ref{eq:bvpmat1} in the new variables is
\begin{align}
\begin{bmatrix}
\dot{Z}_1\\
\dot{Z}_2
\end{bmatrix}
    = U^{-1}P^{-1}NPU\begin{bmatrix}
Z_1\\
Z_2
\end{bmatrix}=\begin{bmatrix}
A_c &B_1\\
0 & -A_c^*
\end{bmatrix}\begin{bmatrix}
Z_1\\
Z_2
\end{bmatrix}.\label{eq:bvpmatZ}
\end{align}
 From the second component of Eqs. \ref{eq:bvpmatZ}, we get $Z_2(t)=e^{-A_c^*t}Z_2(0)$. Since $A_c$ is Hurwitz, $Z_2$ blows up unless we pick $Z_2(0)=0$. This implies $Z_2(t)=0$ for all $t\geq 0$. Using this in the first component of Eqs. \ref{eq:bvpmatZ}, we obtain $Z_1(t)=e^{A_ct}Z_{1}(0)$, and hence, $\lim_{t\rightarrow\infty} Z_1(t)=0$. By invertibility of $PU$, the above two results imply that $\displaystyle\lim_{t\rightarrow\infty}Y_1(t)=0$, and $\displaystyle\lim_{t\rightarrow\infty}Y_2(t)=0$.
\end{proof}
\end{lemma}

Note that if $P=\begin{bmatrix}
P_{11} &P_{12}\\
P_{21} &P_{22}
\end{bmatrix}$, the initial condition $Z_{1}(0)$ above is defined by $(P_{11}+P_{12}X_+)Z_{1}(0)=Y_{10}$. To obtain uniqueness, we need the additional assumption that $(P_{11}+P_{12}X_+)$ is invertible, in which case the unique solution is $(Y_1(t)=(P_{11}+P_{12}X_+)e^{A_ct}(P_{11}+P_{12}X_+)^{-1}Y_{10},Y_2(t)=(P_{21}+P_{22}X_+)e^{A_ct}(P_{11}+P_{12}X_+)^{-1}Y_{10})$. 

From the above results, it follows that if the linearized MFG system of Eq. \ref{eq:mfg_eval} has no eigenvalues on the imaginary axis, and the invertibility assumption holds, the equilibrium solution $(\rho_{\xi},h_{\xi})$ of the MFG system (\ref{eq:mfg_fp},\ref{eq:mfg_hjb}) is linearly stable. 

\subsection{Numerical results}
\begin{figure}[h!]
\subfloat[r=1.4]{
\includegraphics[width=0.24\textwidth]{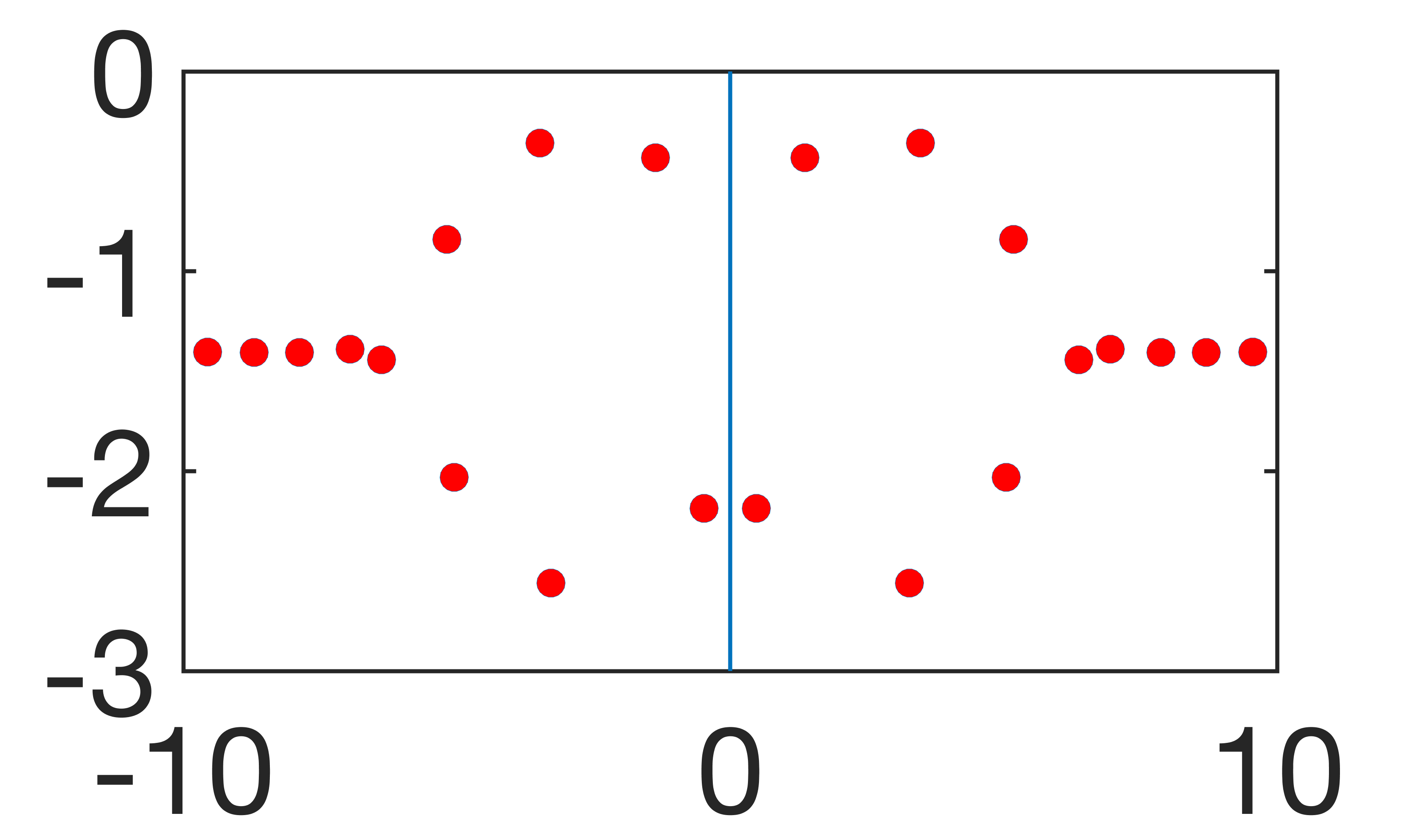}}
\subfloat[r=1.1]{
\includegraphics[width=0.24\textwidth]{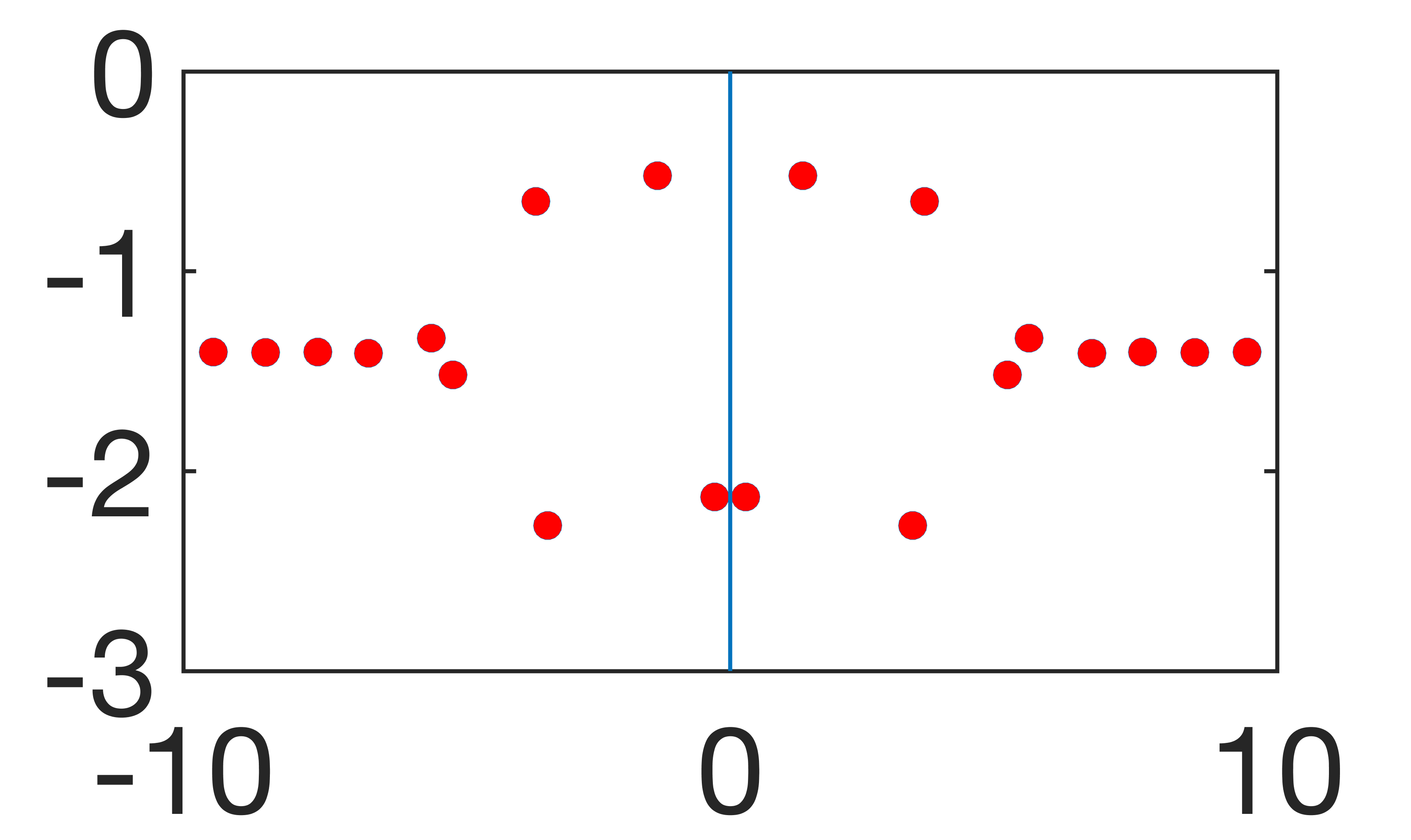}}\\
\subfloat[r=0.95]{
\includegraphics[width=0.24\textwidth]{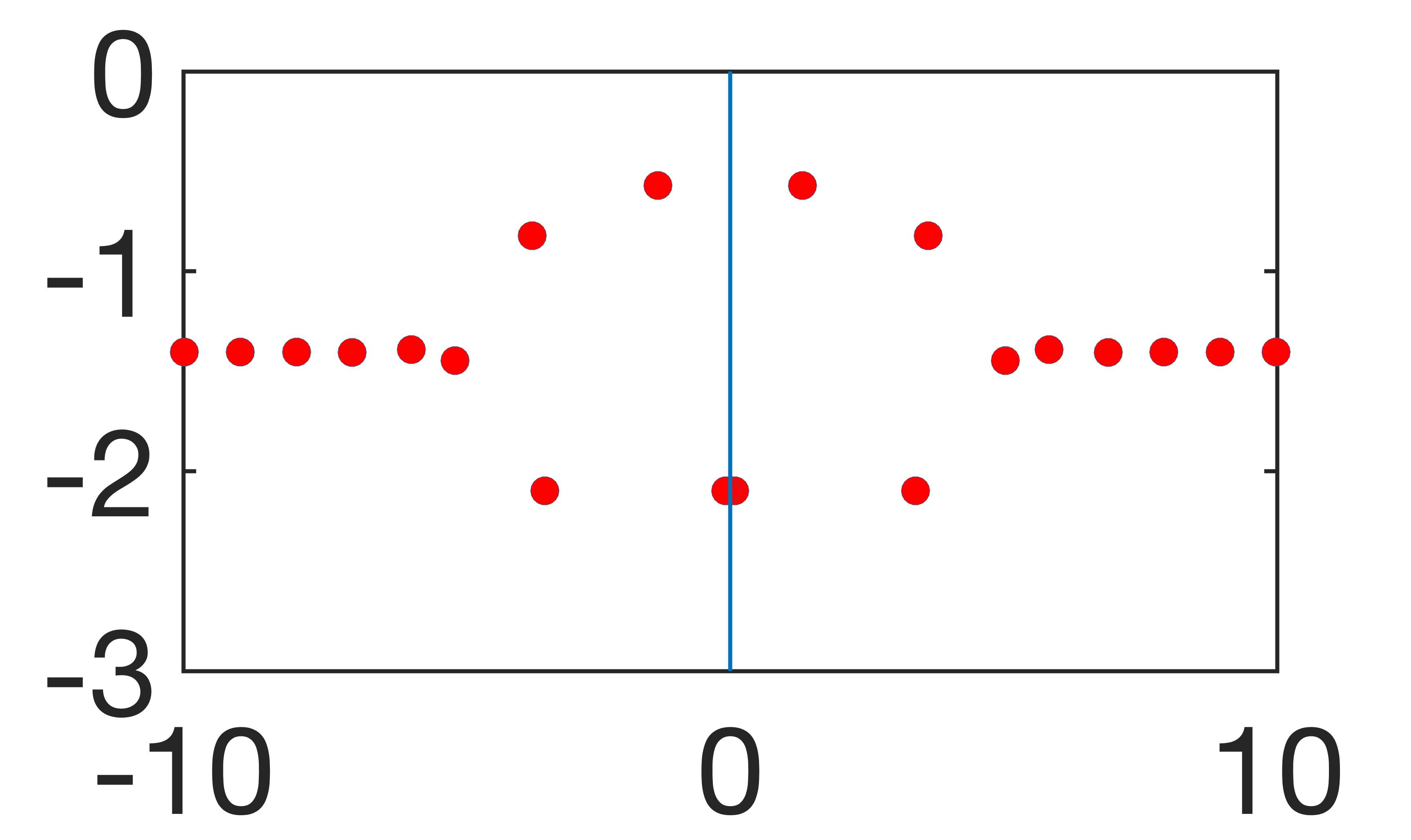}}
\subfloat[r=0.8]{
\includegraphics[width=0.24\textwidth]{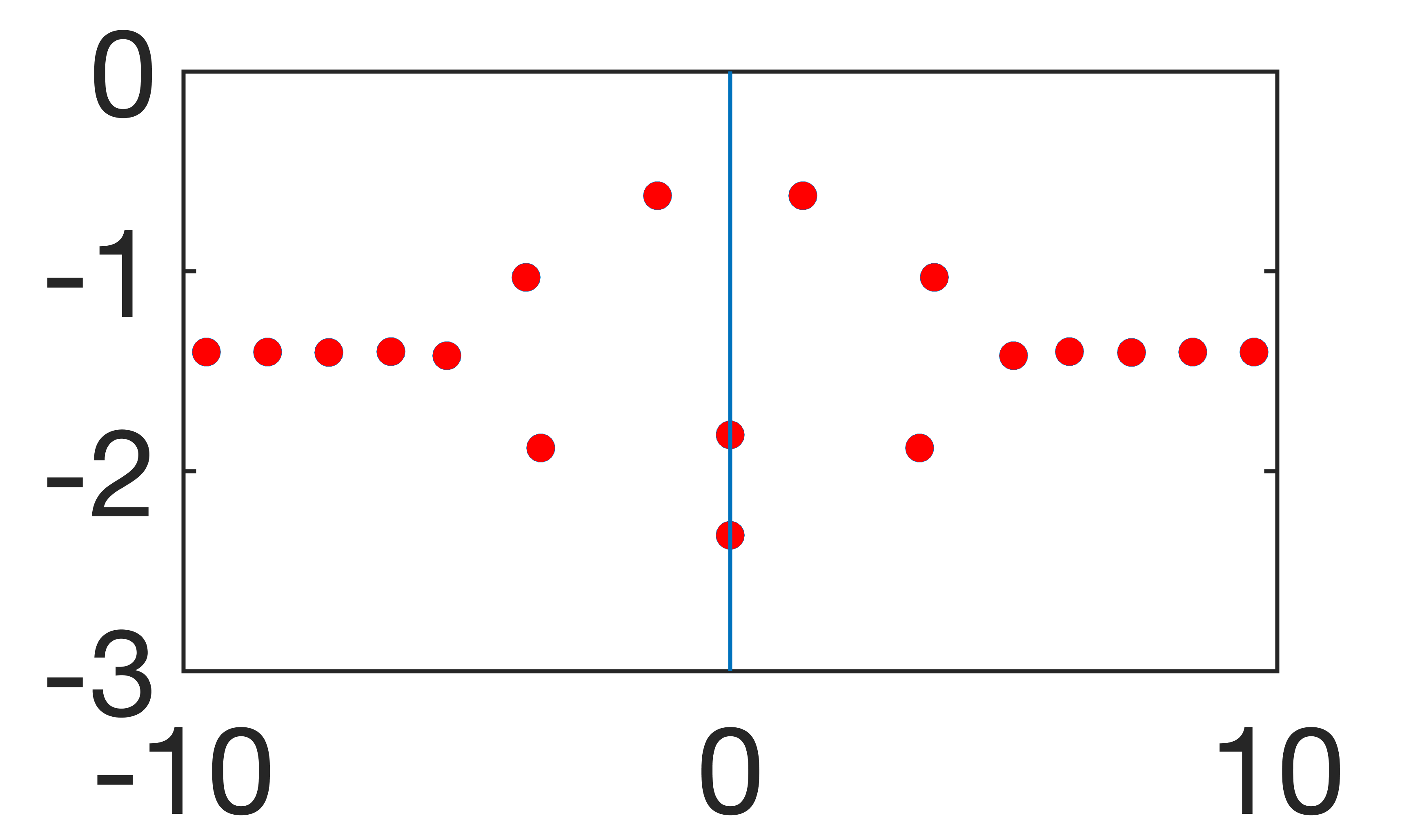}}
\caption{\footnotesize{The spectrum (close to the imaginary axis) of the linearized MFG operator in Eq. \ref{eq:eigsys_mfg} as the control cost $r$ is varied, for Fourier mode $k=1$.  }}
\label{fig:mfg_eval}
\end{figure}

To illustrate our theoretical results, we fix parameters $L=10,h=5,\sigma=2$, and compute the MFG eigenvalues (i.e., eigenvalues of matrix $N$) by numerically solving the algebraic Eq. \ref{eq:mfg_eval} using 22 eigenfunctions for each Fourier mode, i.e., $0\leq p,q\leq 21$. As in the forward equation case, the $k=1$ Fourier mode is the relevant spatial mode for studying stability. Fig. \ref{fig:mfg_eval} shows the eigenspectrum for $k=1$ for various values of the unit control cost $r$. Since the matrix $N$ has the same spectrum as the Hamiltonian matrix $N_s$, this spectrum is symmetric about the imaginary axis. At $r=1.4$, there are no imaginary axis eigenvalues. As $r$ is decreased from 1.4, a pair of eigenvalues approaches the imaginary axis, and eventually collides on it at $r_c\approx 0.95$. As $r$ is reduced futher, the two eigenvalue move away from each other up/down the imaginary axis.

\begin{figure}[h!]
\includegraphics[width=0.45\textwidth]{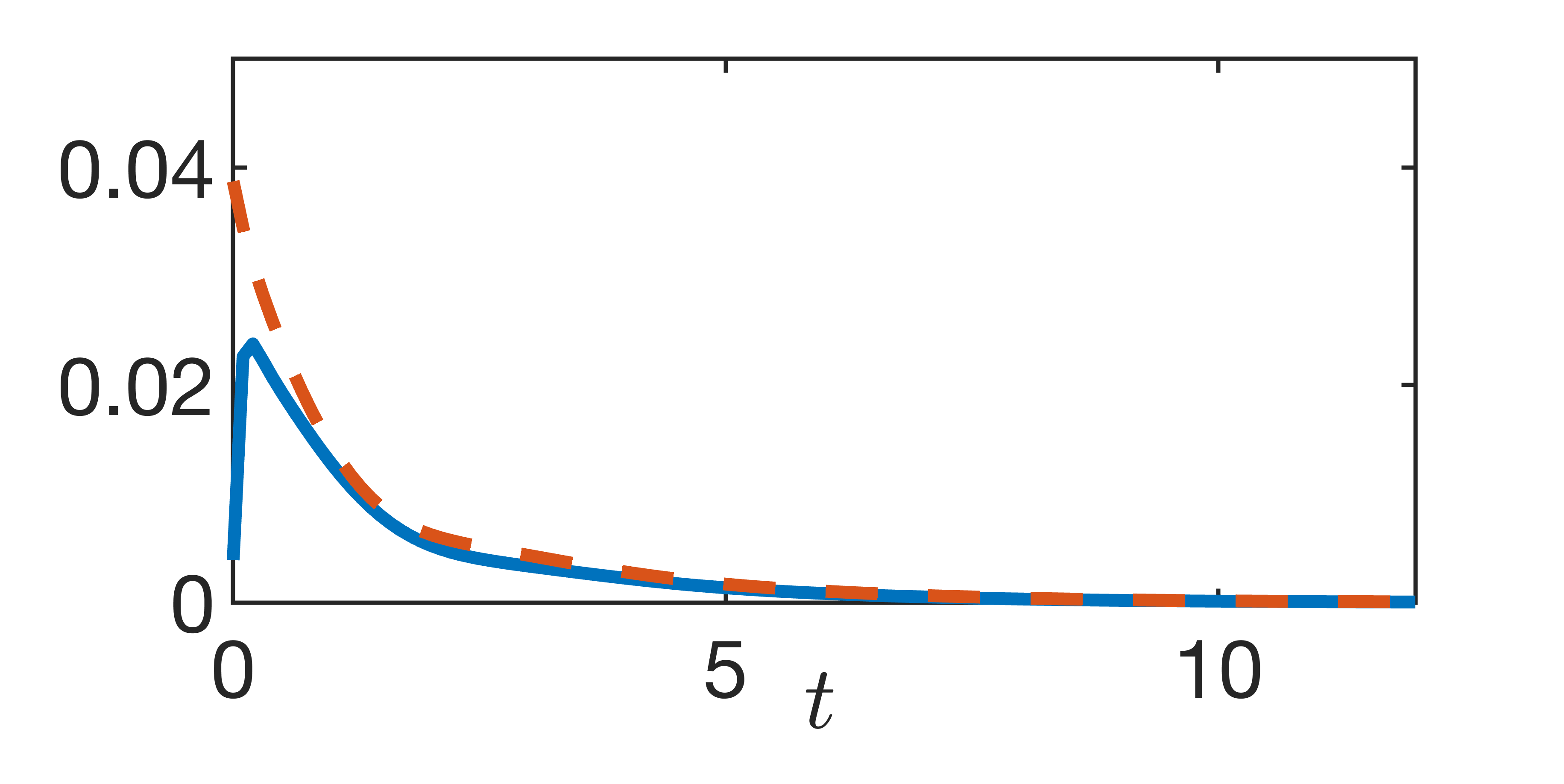}
\caption{\footnotesize{The norm of $Y_1(t)$ (bold) and $Y_2(t)$ (dashed) as a function of time for the unique solution $(Y_1(t),Y_2(t))$ of the BVP of Eq. \ref{eq:bvpmat1}. Here, $r=1.4>r_c$, and $k=1$. We choose an arbitrary $Y_1(0)$, and the corresponding value of $Y_2(0)$ is assigned according to Lemma 1.}}
\label{fig:fb_norm}
\end{figure}

\begin{figure}[h!]
\includegraphics[width=0.45\textwidth]{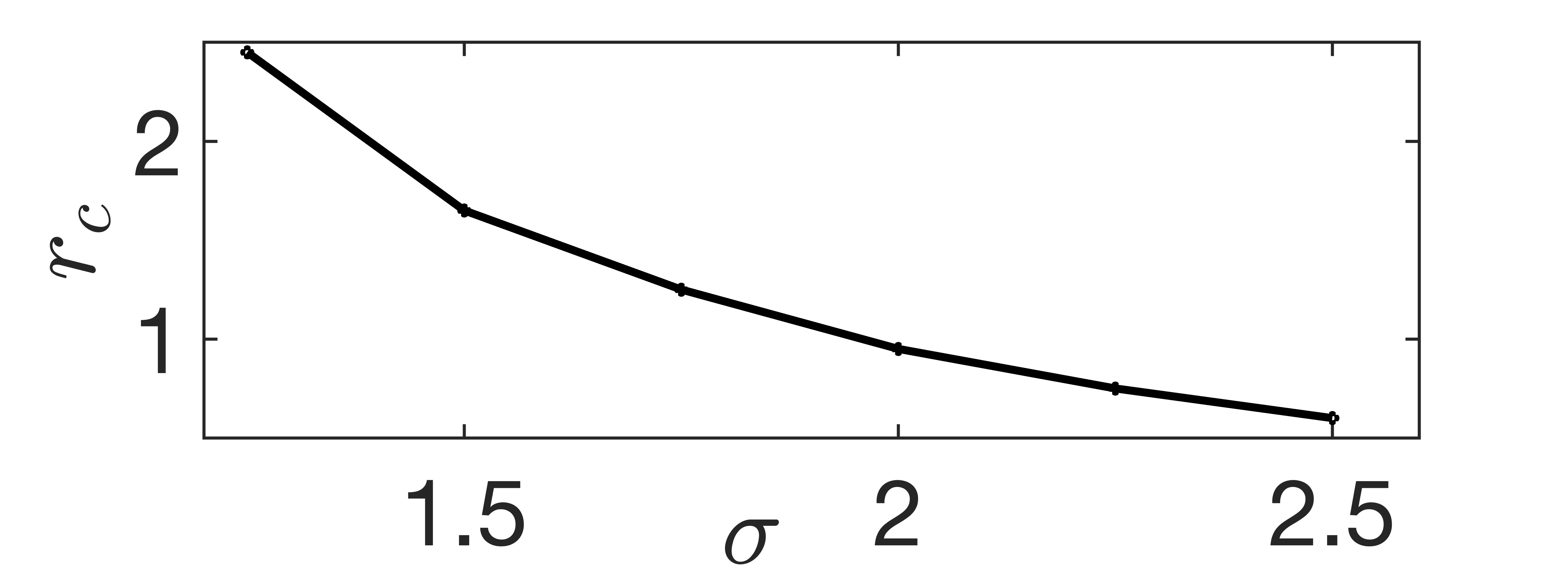}
\caption{\footnotesize{The critical unit control cost $r_c$ as a function of $\sigma$, for $h=5$.}}
\label{fig:sigrc}
\end{figure}

According to Lemma 1, this implies that for $r>r_c$, the MFG equilibrium $(\rho_{\xi},h_{\xi})$ is linearly stable in the sense of Definition 1. We use Schur decomposition to compute $X_+$ using the formula in Theorem 1, and follow the construction of the BVP solutions $(Y_1(t),Y_2(t))$ in Lemma 1. Figure \ref{fig:fb_norm} shows that the norm of the unique solution for an arbitrarily chosen $Y_{10}$ for $r=1.4$ decays to zero. The critical value of unit control cost decreases upon increasing the noise intensity $\sigma$, as shown in Fig. \ref{fig:sigrc}.

\subsubsection*{Travelling wave solutions of the MFG}
\begin{figure}[h!]
\includegraphics[width=0.45\textwidth]{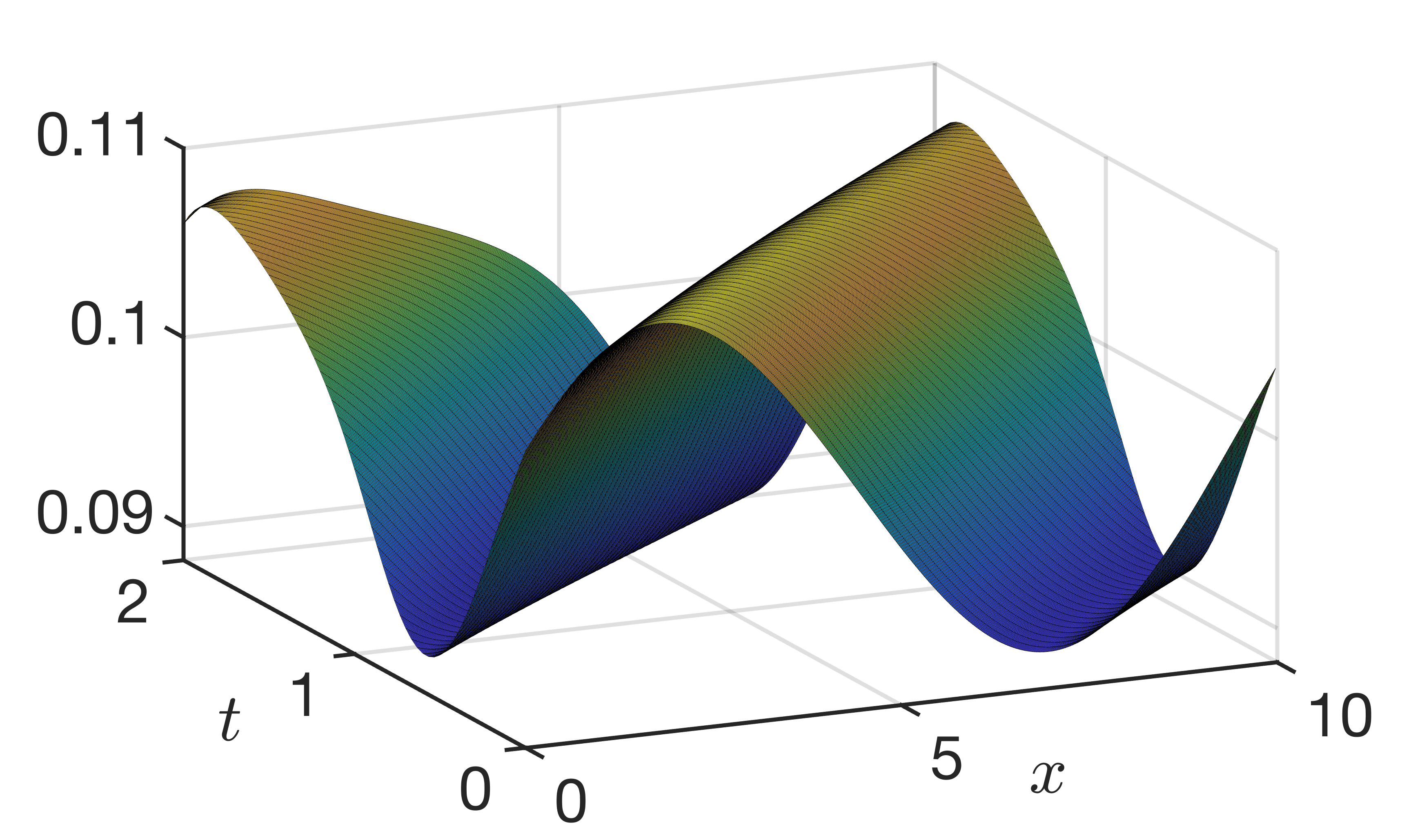}
\caption{\footnotesize{The marginal density $\protect \bigintssss_{\protect \mathbb{R}} \rho(t,x,v)dv$ of a travelling wave solution of the MFG Eqs. (\ref{eq:mfg_fp},\ref{eq:mfg_hjb}) for $r=0.8<r_c$. }}
\label{fig:fb_tw}
\end{figure}

We look for solutions bifurcating from the spatially homogeneous equilibrium $(\rho_{\xi},h_{\xi})$ for $r<r_c$. It is known that MFGs exhibit the turnpike property \cite{zaslavski2005turnpike,cardaliaguet2013long}, i.e., the solution over a large but finite time horizon spends most of its time near the solution of the infinite time problem (the so-called `ergodic' solution). We implment  a Picard-iteration based algorithm for MFGs \cite{lauriere2021numerical} over a large time-horizon in Dedalus, and find that for each $r$ in an open interval $r\in (r_c-\epsilon,r_c)$, the ergodic solution is a travelling wave solution of the form $(\rho(x-\omega t,u),h(x-\omega t,u))$ to the MFG Eqs. (\ref{eq:mfg_fp},\ref{eq:mfg_hjb}). Fig. \ref{fig:fb_tw} shows the marginal density of the travelling wave solution for $r=0.8$. No travelling wave solutions were found for $r>r_c$, in which case the algorithm always converged to one of the two stable equilibria $(\rho_{\pm\xi^*},h_{\pm\xi^*})$.

\section{Conclusions and Discussion}
Via the MFG framework, we have studied transitions between distinct collective behaviors in a population of non-cooperative motile inertial agents that are interacting with their neighbors over a finite distance, and minimizing a biologically inspired cost function. The kinetic MFG model is shown to mimic phase transitions previously observed in the phenomenological Czir\'ok model. 

The linear stability of the equilibrium states is equivalent to the existence of a unique decaying solution to the linearized BVP (in time) derived from the nonlinear MFG PDE system. We provide conditions on the spectrum of the linear operator for such a solution to exist. The explicit calculations are carried out using Fourier-Hermite discretization of the linearized PDE, and use properties of Hamiltonian matrices and Riccati equations. The existence of a traveling wave solution of the MFG when the equilibrium loses stability is shown numerically. A rigorous bifurcation analysis will be taken up in a future work.

While non-equilibrium systems cannot be described by a variational principle, the MFG inverse modelling approach adopted here is based on a generalized optimality principle, which can potentially be extended to kinetic and hydrodynamic descriptions of other systems with decision-making agents. We plan to study some of these extensions in the near future.

\section*{Acknowledgments}
This material is based upon work supported by the US National Science Foundation under Grant No. 2102112.
\section*{Data Availability}
The data that support the findings of this study are available from the corresponding author upon reasonable request.

\clearpage

\bibliography{refs}

\end{document}